\documentclass[final,leqno,onefignum]{siamltex704}

\usepackage{graphicx,epsfig,amsfonts,latexsym,amsmath,enumerate}

\usepackage{pstricks,pst-grad,pst-plot,pst-node}
\usepackage{hyperref}

\RequirePackage{ifthen}
\makeatletter
\newcommand{\logmessage}[1]{\@latex@warning{#1}}
\makeatother
\IfFileExists{../Bibliographie/Useru_Guide.txt}{
  }{
  \IfFileExists{../../Bibliographie/Useru_Guide.txt}{
    }{
    \IfFileExists{../../../Bibliographie/Useru_Guide.txt}{
      }{
        \IfFileExists{../../../../Bibliographie/Useru_Guide.txt}{
          }{
            \IfFileExists{../../../../../Bibliographie/Useru_Guide.txt}{
              }{
        \logmessage{Directory 'Bibliographie' not found}
          }}}}}
\numberwithin{equation}{section}
\newtheorem{remark}[theorem]{Remark}
\newtheorem{assumption}[theorem]{Assumption}
\newtheorem{example}[theorem]{Example}
\newtheorem{algorithm}[theorem]{Algorithm}

\newcommand{\R}{\mathbb{R}}

\newcommand{\half}{\mbox{$\frac{1}{2}$}}

\newcommand{\ignore}[1]{}

\newcommand{\notmid}{\mid\kern-0.5em\not\kern0.5em}

\newcommand{\yd}{y}

\newcommand{\dst}{\displaystyle}

\newcommand{\ba}{\begin{array}}
\newcommand{\ea}{\end{array}}

\title{A convergence analysis of a multi-level projected steepest descent
       iteration for nonlinear inverse problems in Banach spaces subject
       to stability constraints
   \thanks{This research was supported in part National Science
     Foundation grant CMG DMS-1025318 and in part by
     the members of the Geo-Mathematical Imaging Group at Purdue
     University. The work of OS has been supported by the Austrian
     Science Fund (FWF) within the national research networks
     Photo\-acoustic Imaging in Biology and Medicine, project
     S10505 and Geometry and Simulation S11704.}
               }

\author{Maarten V. de Hoop\thanks{Center for Computational and Applied
    Mathemematics, Purdue University, West Lafayette, IN 47907 ({\tt
      mdehoop@purdue.edu}).}
\and Lingyun Qiu\thanks{Center for Computational and Applied
  Mathemematics, Purdue University, West Lafayette, IN 47907 ({\tt
    qiu@purdue.edu}).}
\and Otmar Scherzer\thanks{Computational Science Center, University of
  Vienna, Nordbergstr. 15, A-1090 Vienna, Austria ({\tt
    otmar.scherzer@univie.ac.at}).}  }

\begin{document}

\maketitle

\setcounter{page}{1}

\begin{abstract}
We consider nonlinear inverse problems described by operator equations
in Banach spaces. Assuming conditional stability of the inverse
problem, that is, assuming that stability holds on a closed, convex
subset of the domain of the operator, we introduce a novel nonlinear
projected steepest descent iteration and analyze its convergence to an
approximate solution given limited accuracy data. We proceed with
developing a multi-level algorithm based on a nested family of closed,
convex subsets on which stability holds and the stability constants
are ordered. Growth of the stability constants is coupled to the increase in
accuracy of approximation between neighboring levels to ensure that the
algorithm can continue from level to level until the iterate satisfies
a desired discrepancy criterion, after a finite number of steps.
\end{abstract}


\thispagestyle{plain}
\markboth{M.V. DE HOOP, L. QIU AND O. SCHERZER}{MULTI-LEVEL NONLINEAR
  PROJECTED STEEPEST DESCENT ITERATION}

\section{Introduction}
\label{sec:intro}

We consider nonlinear inverse problems described by operator equations
in Banach spaces. Assuming conditional stability of the inverse
problem, we introduce a nonlinear projected steepest descent iteration
and analyze its convergence. We take the point of view of
reconstructing an approximation of the solution to the inverse problem
in a closed, convex subset of the domain on which the operator is
defined and where the stability holds. Assuming that we can identify
a nested sequence of closed, convex subsets on which the stability
holds such that the stability constant grows in a controlled way, we then extend our analysis to a multi-level approach which
mitigates this growth via successive approximation. We account also for the possibility that a parameter in the operator which defines the inverse problem, and changes the data, affects for a given closed, convex subset the accuracy of approximation, and the stability constant, as well. Our multilevel approach results in a radius of convergence which is significantly larger than the one in the single level approach. In
our analysis, we incorporate inaccuracy of the data. Our analysis
applies, for example, to electrical impedance tomography (EIT) and
inverse boundary value problems for the Helmholtz equation using
multiple frequencies.

Initially, we consider a class of inverse problems defined by a
nonlinear map from parameter or model functions to the data. The
parameter functions and data are contained in certain Banach
spaces. This situation can be modeled mathematically by the operator
equation
\begin{equation}\label{forward-equ}
   F(x) = y , \quad x \in \mathcal{D}(F) ,\ y \in Y ,
\end{equation}
with domain $\mathcal{D}(F) \subset X$, where $X$ and $Y$ are Banach
spaces.
We assume that $F$ is continuous, and that $F$ is locally Fr\'{e}chet
differentiable. We do \emph{not} assume that the data are
attainable, that is, $\yd$ may not belong to the range of $F$. We
assume that there exists a closed, convex subset $Z \subset X$ such
that
\begin{equation}\label{eq:constab}
   \Delta_p(x,\tilde{x}) \le \mathfrak{C}^p \| F(x) - F(\tilde{x}) \|^p ,\quad
             \forall x, \tilde{x} \in Z.
\end{equation}
Here $\Delta_p$ denotes the Bregman distance (defined below) and $p>1$. This states conditional Lipschitz stability of the
inverse problem. Motivated by \cite{Hoop2012}, we employ a steepest
descent iteration, here, to give an approximation to the solution of
(\ref{forward-equ}). More precisely, we construct a sequence of
parameter functions by a projected gradient descent iteration with
posterior stepsize.

In many inverse problems, logarithmic type stability is the optimal
stability obtained with minimal assumptions on the domain or pre-image
space; see, for example, \cite{Mandache2001}. By constraining the
pre-image space, however, Lipschitz stability can be obtained; for
the case of EIT, see \cite{Alessandrini2005,Beretta2011} and for the
case of inverse boundary value problems for the Helmholtz equation,
see \cite{Beretta2012,Ammari2012}. This is reflected by conditional stability given
in (\ref{eq:constab}). The mentioned projected gradient descent
iteration can then be viewed as a projection regularization method,
which is natural and avoids possibly artificial regularization
techniques \cite{Kaltenbacher2006}.

Our first main result concerns restricted convergence of the
projected steepest descent iteration with a certain
Lipschitz type stability condition on a closed, convex
subset. Moreover, we prove monotonicity of the residuals defined by
the sequence induced by the iteration.  This result is related to two
areas of iterative regularization, which are steepest descent
algorithms for solving nonlinear inverse problems
\cite{Neubauer1995,Scherzer1996,Kaltenbacher2008} and projected iteration
regularization techniques for the solution of inverse problems with
convexity constraints. The latter have been analyzed mostly in the
context of \emph{linear} inverse problems (see, for example,
\cite{Eicke1992}) and later as accelerated methods in
\cite{Daubechies2008}. Accelerated methods have been modified to nonlinear
problems by \cite{Teschke2010}. The main differences of our work to
the above mentioned papers are the conditions under which we prove
convergence. In fact, instead of source and nonlinearity conditions
(as in \cite{Neubauer1995,Scherzer1996}), we assume certain H\"older or Lipschitz
stability of the inverse problem. This is a novel view point, which
has been raised in \cite{Hoop2012}. The steepest descent method proposed
here is a generalization of the steepest descent method for unconstrained linear
problems (see for example \cite{Gilyazov1977}). It is however different from the generalization for nonlinear problems proposed
in \cite{Neubauer1995,Scherzer1996}, even for unconstrained problems.

Based on our first main result, we then introduce a multilevel
algorithm. We assume that there are closed, convex subsets $\{
Z_{\alpha} \}_{\alpha \in \R}$ of $X$, on which the restricted
operator $F_{\alpha} = F \mid_{Z_{\alpha}}$ exhibits a certain
H\"{o}lder or Lipschitz type stability estimate with stability
constant $\mathfrak{C}_{\alpha}$, that is,
\begin{equation}\label{eq:constab_alpha}
   \Delta_p(x,\tilde{x}) \le \mathfrak{C}_{\alpha}^p
     \| F_{\alpha}(x) - F_{\alpha}(\tilde{x}) \|^p ,\quad
              \forall x, \tilde{x} \in Z_{\alpha} .
\end{equation}
In fact, $F_{\alpha}$ need not be a restriction of $F$ only, but can
also account for a varying parameter in $F$ which does affect the
data. Here, we assume that $Z_{\alpha_1} \subset Z_{\alpha_2}$ and
$\mathfrak{C}_{\alpha_1} \le \mathfrak{C}_{\alpha_2}$ if $\alpha_1 < \alpha_2$. In the
context of discretization methods, $Z_{\alpha}$ stands for a
finite-dimensional subspace of $X$ and the number of basis vectors
increases as $\alpha$ increases, while the projection can be an
orthogonal projection on $Z_\alpha$. In our second main result, we
introduce a condition on the stability constants and on the
approximation errors between neighboring levels. These conditions between levels are coupled and
guarantee that the result from the previous level is a proper
starting point for the present level. Thus, the algorithm can continue
from level to level until the desired discrepancy criterion is
satisfied.

\medskip\medskip
\section{Preliminaries}
\label{sec:Preliminaries}
Several constants appear in the analysis. For the readers convenience we have grouped them as follows:
\begin{enumerate}
\item $\mathfrak{C}$ denotes a constant for the Lipschitz  stability of the inverse mapping of $F$ (cf. (\ref{eq:constab}), (\ref{eq:constab_alpha})), \item $\mathfrak{L}$ and $\hat{\mathfrak{L}}$ are properties of the the operator $F$ (cf. (\ref{eq:Lipschitz}) and (\ref{eq:Lipschitz_derivative})).
\item $C$ and $G$ with and without subscripts denote properties of the Banach space (cf. (\ref{Bregman-norm-rela1}), (\ref{Bregman-norm-rela2})).
\end{enumerate}

\subsection{Duality mappings}
\label{subsec:Duality}

Let $X$ and $Y$ be Banach spaces. The duals of $X$ and $Y$ are denoted
by $X^*$ and $Y^*$, respectively. Their norms are denoted uniformly by
$\|\cdot\|$. We denote the space of continuous linear operators $X
\rightarrow Y$ by $\mathcal{L}(X,Y)$. Let $F : \mathcal{D}(F) \subset
X \rightarrow Y$ be continuous. Here $\mathcal{D}(F)$ denotes the
domain of definition of the nonlinear operator $F$. Let $h \in \mathcal{D}(F)$ and $k \in X$ and assume
that $h + t(k-h) \in \mathcal{D}(F)$ for all $t\in (0, t_0)$ for some
$t_0>0$, then we denote by $DF(h)(k)$ the directional derivative of $F$
at $h\in \mathcal{D}(F)$ in direction $k\in \mathcal{D}(F)$, that is,
\[
DF(h)(k):= \lim_{t \rightarrow 0^+} \frac{F(h+tk) - F(h)}{t}.
\]
If $DF(h) \in \mathcal{L}(X,Y)$, then $F$ is called G\^{a}teaux
differentiable at $h$. If, in addition, the limit is uniform for all $k$ belonging a neighborhood of $0$, $F$ is called Fr\'{e}het
differentiable at $h$.
For $x \in X$ and $x^* \in X^*$, we write the dual pair as $\langle x,x^* \rangle
= x^*(x)$. For a linear operator $A\in \mathcal{L}(X,Y)$, we write $A^*$ for the dual operator $A^* \in
\mathcal{L}(Y^*,X^*)$ and $\|A\| = \|A^*\|$ for the operator norm of
$F$. We let $1 < p,q < \infty$ be conjugate exponents, that is,
\begin{equation*}
   \frac{1}{p} + \frac{1}{q} = 1 .
\end{equation*}
For $p > 1$, the subdifferential mapping $J_p = \partial f_p : X
\rightarrow 2^{X^*}$ of the convex functional $f_p: x\mapsto
\frac{1}{p}\|x\|^p$ defined by
\begin{equation}\label{Duality-map-def}
   J_p(x) = \{ x^* \in X^* \mid \langle x, x^* \rangle
           = \|x\| \cdot \|x^*\| \mbox{ and } \|x^*\| = \|x\|^{p-1}\}
\end{equation}
is called the duality mapping of $X$ with gauge function $t \mapsto t^{p-1}$. Generally, the duality mapping is set-valued. In order to let $J_p$ be single valued, we need to introduce the notion of convexity and smoothness of Banach spaces.

One defines the convexity modulus $\delta_X$ of $X$ by
\begin{equation}\label{def:convex-module}
   \delta_X(\epsilon) = \inf_{ x, \tilde{x} \in X} \{
      1 - \| \half (x + \tilde{x})\| \mid \|x\| = \|\tilde{x}\| = 1
            \mbox{ and }\|x - \tilde{x}\| \ge \epsilon\}
\end{equation}
and the smoothness modulus $\rho_X$ of $X$ by
\begin{equation}\label{def:smooth-module}
   \rho_X(\tau) = \sup_{x, \tilde{x} \in X} \{
      \half (\|x + \tau \tilde{x}\| + \|x - \tau \tilde{x}\|
           - 2) \mid \|x\| = \|\tilde{x}\| = 1 \} .
\end{equation}

\medskip\medskip

\begin{definition}
A Banach space $X$ is said to be
\begin{enumerate}[(a)]
\item uniformly convex if there exists an $\epsilon \in (0,2]$ such that $\delta_X(\epsilon) > 0$\,,
\item uniformly smooth if $\lim_{\tau \rightarrow 0} \frac{\rho_X(\tau)}{\tau} =0$,
\item convex of power type $p$ or $p$-convex if there exists a constant $C > 0$ such that $\delta_X(\epsilon) \ge C \epsilon^p$,
\item smooth of power type $q$ or $q$-smooth if there exists a constant $C > 0$ such that $\rho_X(\tau) \le C \tau^q$.
\end{enumerate}
\end{definition}
For a detailed introduction to the geometry of Banach spaces and the duality mapping, we refer to \cite{Cioranescu1990,Schopfer2006}. We list the properties we need here in the following theorem.

\medskip\medskip

\begin{theorem} Let $p>1$. The following statements hold true:
\begin{enumerate}[(a)]
\item For every $x \in X$, the set $J_p(x)$ is not empty and it is convex and weakly closed in $X^*$.
\item Theorem of Milman-Pettis: If a Banach space is uniformly convex, it is reflexive.
\item A Banach space $X$ is uniformly convex (resp. uniformly smooth) if and only if $X^*$ is uniformly smooth (resp. uniformly convex).
\item If a Banach space $X$ is uniformly smooth, $J_p(x)$ is single valued for all $x \in X$.
\item If a Banach space $X$ is uniformly smooth and uniformly convex, $J_p(x)$ is bijective and the inverse $J_p^{-1}: X^* \rightarrow X$ is given by $J_p^{-1} = J_q^*$ with $J_q^*$ being the duality mapping of $X^*$ with gauge function $t \mapsto t^{q-1}$, where $1 < p,q < \infty$ are conjugate exponents.
\end{enumerate}
\end{theorem}

\medskip\medskip

\subsection{Bregman distances}
\label{subsec:Bregman}

Because the geometrical characteristics of Banach spaces are different
from those of Hilbert spaces, it is often more appropriate to use the
Bregman distance instead of the conventional norm-based functionals
$\|x-\tilde{x}\|^p$ or $\|J_p(x) - J_p(\tilde{x})\|^p$ for convergence
analysis. This idea goes back to Bregman
\cite{Bregman1967}. 

\begin{definition}
Let $X$ be a uniformly smooth Banach space and $p > 1$. The Bregman
distance $\Delta_p (x, \cdot)$ of the convex functional $x\mapsto
\frac{1}{p}\|x\|^p$ at $x \in X$ is defined as
\begin{equation}\label{def:Bregman-distance}
   \Delta_p(x,\tilde{x}) =
       \frac{1}{p} \|\tilde{x}\|^p - \frac{1}{p} \|x\|^p
    - \langle J_p(x), \tilde{x} - x\rangle , \quad \tilde{x} \in X ,
\end{equation}
where $J_p$ denotes the duality mapping of $X$ with gauge function $t \mapsto t^{p-1}$. Note, that under the general assumptions of this paper the duality mapping $J_p$ is single valued.
\end{definition}

\medskip\medskip

In the following theorem, we summarize some facts concerning the
Bregman distance and the relationship between the Bregman distance and the
norm \cite{Alber1997,Alber1996,Butnariu2000,Xu1991}.

\medskip\medskip

\begin{theorem}\label{thm:Bregman-distance}
Let $X$ be a uniformly smooth and uniformly convex Banach space. Then,
for all $x, \tilde{x} \in X$, the following holds:
\begin{enumerate}[(a)]
\item
\begin{eqnarray}\label{eq:BDpq}
   \Delta_p(x,\tilde{x}) &=& \frac{1}{p}\|\tilde{x}\|^p
   - \frac{1}{p}\|x\|^p - \langle J_p(x), \tilde{x} \rangle + \|x\|^p
\\
   &=& \frac{1}{p}\|\tilde{x}\|^p + \frac{1}{q}\|x\|^p
       - \langle J_p(x), \tilde{x}\rangle .
\nonumber
\end{eqnarray}
\item $\Delta_p(x,\tilde{x}) \ge 0$ and $\Delta_p(x,\tilde{x}) = 0
  \Leftrightarrow x = \tilde{x}.$
\item $\Delta_p$ is continuous in both arguments.
\item The following statements are equivalent
\begin{enumerate}[(i)]
\item $\lim_{n \rightarrow \infty} \|x_n - x\| = 0, $
\item $\lim_{n \rightarrow \infty} \Delta_p(x_n , x) = 0 ,$
\item $\lim_{n \rightarrow \infty} \|x_n\| = \|x\|$ and $ \lim_{n
  \rightarrow \infty} \langle J_p(x_n), x\rangle = \langle J_p(x),
  x\rangle$.
\end{enumerate}
\item If $X$ is $p$-convex, there exists a constant $C_p > 0$ such that
\begin{equation}\label{Bregman-norm-rela1}
   \Delta_p(x,\tilde{x}) \ge \frac{C_p}{p} \|x-\tilde{x}\|^p .
\end{equation}
\item If $X^*$ is $q$-smooth, there exists a constant $G_q > 0$ such that
\begin{equation}\label{Bregman-norm-rela2}
   \Delta_q(x^*,\tilde{x}^*) \le \frac{G_q}{q} \|x^*-\tilde{x}^*\|^q ,
\end{equation}
for all $x^*, \tilde{x}^*\in X^*.$
\end{enumerate}
\end{theorem}

\medskip\medskip

The Bregman distance $\Delta_p$ is similar to a metric, but, in
general, does not satisfy the triangle inequality nor symmetry. In a
Hilbert space, $\Delta_2(x,\tilde{x}) = \frac{1}{2}
\|x-\tilde{x}\|^2.$

\subsection{Bregman Projection}
\label{subsec:Projection}

In this subsection, we briefly introduce the Bregman projection and
its properties, especially, the non-expansiveness. A comprehensive introduction to this topic, including a proof of Lemma~\ref{lemma:proj}, can be found in \cite{Butnariu2000}.

\begin{definition}\label{def:projection}
Let $X$ be a uniformly smooth Banach space and $p > 1$. Given a closed
convex set $Z\subset X$ and Bregman distance $\Delta_p$, which is
defined in Definition~\ref{def:Bregman-distance}, the Bregman
projection of a point $x\in X$ onto $Z$ is the point
\begin{equation}\label{projection}
P_Z(x) = \arg\min \{ \Delta_p(y,x) \mid y\in Z\}.
\end{equation}
\end{definition}


\medskip\medskip

\begin{definition}\label{def:t-nonexpan}
Let $T:X \rightarrow X$ be an operator. The point $z\in X$ is called a non-expansivity pole of $T$ if, for every $x\in X$,
\[
\Delta_p(T(x),T(z)) + \Delta_p(x,T(x)) \le \Delta_p(x,z).
\]
A operator $T$, which has at least one non-expansivity pole, is called totally non-expansive.
\end{definition}

\begin{lemma}\label{lemma:proj}
Let $X$ be a uniformly smooth Banach space and $p > 1$ and $Z\subset X$ be a closed convex subset. The following statements hold:
\begin{enumerate}[(a)]
\item The Bregman projection $P_Z$ is well defined;
\item $P_Z$ is totally non-expansive and every point in $Z$ is a non-expansivity pole of $P_Z$;
\item For every $z\in Z$,
\begin{equation}\label{non-expansive}
\Delta_p(P_Z(x), z) \le \Delta_p( x , z), \quad \forall x\in X.
\end{equation}
\end{enumerate}
\end{lemma}

\medskip\medskip

Throughout this paper, we assume that $X$ is $p$-convex and $q$-smooth
with $p,q > 1$, and hence it is uniformly smooth and uniformly
convex. Furthermore, $X$ is reflexive and its dual $X^*$ has the same
properties. $Y$ is allowed to be an arbitrary Banach space; $j_p$ will
be a single-valued selection of the possibly set-valued duality
mapping of $Y$ with gauge function $t \mapsto t^{p-1}$, $p >
1$. Further restrictions on $X$ and $Y$ will be indicated in the
respective theorems below.

\section{Convergence of a projected steepest descent iteration}
\label{sec:convergence}

Here, we assume conditional stability, that is stability if operator
$F$ is restricted to a closed, convex subset, $Z$, of $X$(see (\ref{stab-Banach})). We introduce
a projected steepest descent iteration and analyze its convergence. In
this section, we keep $Z$ fixed. We are concerned with an approximate
solution, in $Z$, of the inverse problem subject to a discrepancy
principle.

\medskip\medskip

\begin{assumption}\label{assumption:forward-operator}
Let
\[
\mathcal{B}=\mathcal{B}^{\Delta}_{\rho}(z^{\dagger}) = \{ x \in
X\ |\ \Delta_p(x,z^{\dagger}) \le \rho \} \subset \mathcal{D}(F)
 \]
 for some $\rho>0$, where $\rho$ here will come into play as a convergence radius and $z^\dagger$ is defined below.
\begin{enumerate}[(a)]
\item The Fr\'{e}chet derative, $DF$, of $F$ is Lipschitz
  continuous on $\mathcal{B}$ and
    \begin{equation}\label{eq:Lipschitz}
   \| DF(x) \| \le \hat{\mathfrak{L}}
   \quad \forall x\in \mathcal{B},
\end{equation}
  \begin{equation}\label{eq:Lipschitz_derivative}
   \| DF(x) - DF(\tilde{x}) \| \le \mathfrak{L} \| x - \tilde{x} \|
   \quad \forall x, \tilde{x} \in \mathcal{B}.
\end{equation}
\item $F$ is weakly sequentially closed, i.e.,
\begin{equation*}
   \left. \begin{array}{rl}
      & x_n \rightharpoonup x ,
\\
      & F(x_n) \rightarrow y
   \end{array} \right\} \Rightarrow
   \left\{ \begin{array}{rl}
      & x \in \mathcal{D}(F) ,
\\
      & F(x) = y .
   \end{array} \right.
\end{equation*}
\item Let $Z$ denote a closed, convex subset of $X$. The inversion has
  the uniform Lipschitz type stability for elements in $Z$, i.e.,
  there exists a constant $\mathfrak{C} > 0$ such that
\begin{equation}\label{stab-Banach}
   \Delta_p(x,\tilde{x}) \le \mathfrak{C}^p \|F(x) - F(\tilde{x})\|^p
   \quad \forall x, \tilde{x} \in \mathcal{B} \cap Z.
\end{equation}
\end{enumerate}
\end{assumption}

\medskip\medskip

For given data $\yd \in Y$, we assume that
\begin{equation}\label{eta}
   \operatorname{dist}(\yd , F(Z)) \le \eta,
\end{equation}
for some $\eta >0$. Note that $F$ is continuous and $Z$ is closed. Hence there must exist a $z^\dag \in Z$ such that
\begin{equation}
   \| F(z^\dag) - \yd \| = \operatorname{dist}(\yd , F(Z)) .
\end{equation}
Note, that this condition also accounts for data errors.

The stopping index $K = K(\eta)$
of the following iteration is determined by a discrepancy principle
\begin{equation}\label{discrepancy-principle}
   K(\eta) := \min \{ k \in\mathbb{N}  \mid \|F(x_k) - y\| \le \hat{\eta} \}
\end{equation}
with a fixed
\begin{equation}
\label{db}
\hat{\eta} > 3 \eta\;.
\end{equation}
We introduce the following algorithm:
\medskip
\begin{algorithm}\label{algo:1}
We fix some abbreviations first: For $x_k$, $k=0,1,2,\ldots$, fixed denote
\begin{equation}
\label{abbreviation}
     R_k = F(x_k) - \yd\,,\quad T_k  = DF(x_k)^*j_p(F(x_k) - \yd)\,, \quad
     r_k = \|R_k\|\,, \quad t_k = \|T_k\|\;.
\end{equation}
Moreover, we define
\begin{equation}
\label{abbreviation_2a}
     \tilde{\mathfrak{C}} \, := \frac{1}{2} \left(\frac{C_p}{p}\right)^{-2/p}\mathfrak{L}\mathfrak{C}^2\,,
\end{equation}
and for $k=0,1,\ldots$
\begin{equation}
\label{abbreviation_2}
\begin{aligned}
     \hat{t}_k & := G_q t_k^q\,,\\
     u_k & := -\tilde{\mathfrak{C}} r_k^2 + (1 - 2\tilde{\mathfrak{C}}\eta)r_k -\eta -\tilde{\mathfrak{C}} \eta^2 \,,\\
     v_k & := \hat{t}_k^{-\frac{1}{q-1}} u_k^{\frac{1}{q-1}}r_k^{p^2-p}(r_k - \eta) - \frac{1}{q} \hat{t}_k^{-\frac{1}{q-1}}u_k^pr_k^{p^2-p} \,,\\
     w_k & :=\frac{\mathfrak{L}}{2}\left(\frac{C_p}{p}\right)^{-2/p}  \hat{t}_k^{-\frac{1}{q-1}}u_k^{\frac{1}{q-1}} r_k^{p^2-p} \,,\\
     \mu_k & := \hat{t}_k^{-\frac{1}{q-1}} u_k^{\frac{1}{q-1}} r_k^{\frac{p-1}{q-1}} \;.
\end{aligned}
\end{equation}

Now, the main steps of the algorithm:
\begin{enumerate}[$(S1)$]
\setcounter{enumi}{-1}
\item Choose a starting point $x_0\in Z$ such that
\begin{equation}\label{converge-radius}
\Delta_p(x_0, z^\dag) < \rho := \frac{C_p}{p}
(2\tilde{\mathfrak{C}}\hat{\mathfrak{L}})^{-p} \left(1+\sqrt{1-8\tilde{\mathfrak{C}}\eta} - 4\eta\tilde{\mathfrak{C}}\right)^p ,
\end{equation}
where $z^\dag$ is specified in Theorem~\ref{thm:stopping-rule} below.
\item Compute the new iterate via
\begin{equation}\label{Steepest descent}
\begin{array}{rl}
   \tilde{x}_{k+1}  = & J_q^*(J_p(x_{k}) - \mu_k T_k)
\\[0.2cm]
   x_{k+1}  = & \mathcal{P}_Z(\tilde{x}_{k+1} ) .
\end{array}
\end{equation}
Set $k \leftarrow k + 1$ and repeat step $(S1)$.
\end{enumerate}
\end{algorithm}

\medskip\medskip

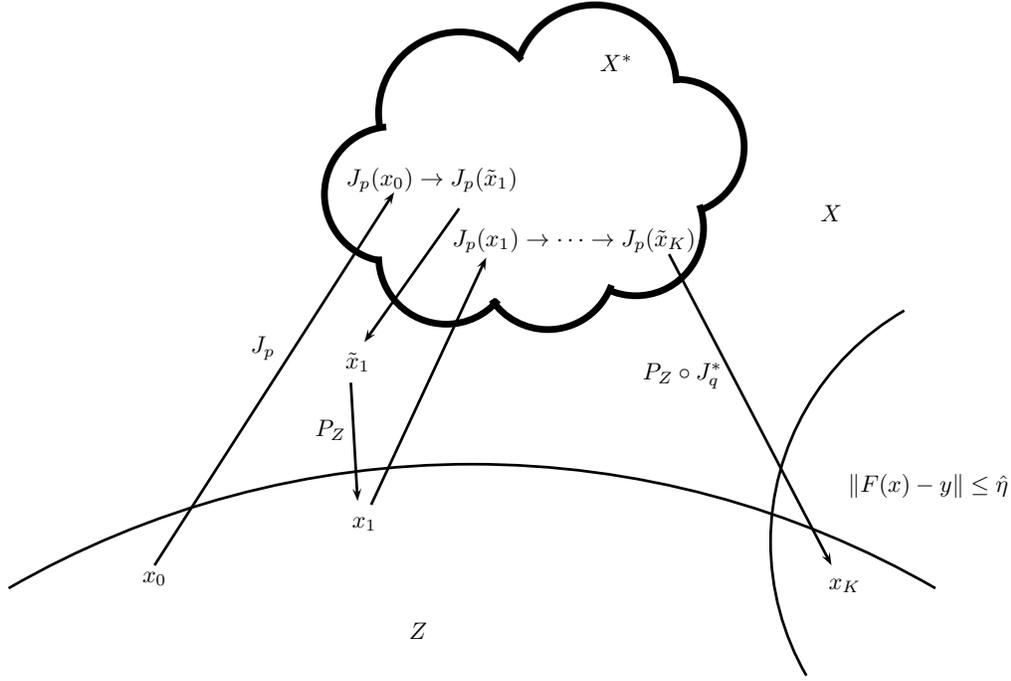
\begin{figure}[htp]
\centering
\scalebox{0.9} 
{
\begin{pspicture}(7.2,7)(21,17)
\psarc[linewidth=0.1](12.510087,14.285){1.0}{95.0}{260.0}
\psarc[linewidth=0.1](13.310087,13.365){1.0}{180.0}{320.0}
\psarc[linewidth=0.1](14.810087,13.285){1.0}{214.0}{340.0}
\psarc[linewidth=0.1](16.110086,13.785){1.0}{245.0}{20.0}
\psarc[linewidth=0.1](16.710087,14.985){1.0}{290.0}{90.0}
\psarc[linewidth=0.1](15.510087,15.885){1.2}{2.0}{160.0}
\psarc[linewidth=0.1](13.510087,15.485){1.2}{40.5}{190.0}
\usefont{T1}{ptm}{m}{n}
\rput(15.822344,16.225555){$ X^*$}
\usefont{T1}{ptm}{m}{n}
\rput(19,14){$ X$}
\usefont{T1}{ptm}{m}{n}
\rput(9.0,8.6){$ x_0$}
\usefont{T1}{ptm}{m}{n}
\rput(10.6,12.005555){$J_p$}
\psline[linewidth=0.04cm,arrowsize=0.05291667cm 2.0,arrowlength=1.4,arrowinset=0.4]{->}(9.0,8.8)(12.539531,14.315556)
\usefont{T1}{ptm}{m}{n}
\rput(13.1,14.5){$J_p(x_0)\rightarrow J_p(\tilde{x}_1)$}
\psline[linewidth=0.04cm,arrowsize=0.05291667cm 2.0,arrowlength=1.4,arrowinset=0.4]{<-}(19,8.8)(16.6,13.4)
\usefont{T1}{ptm}{m}{n}
\rput(19.2,8.5){$x_K$}
\usefont{T1}{ptm}{m}{n}
\psarc[linewidth=0.04](13.69,-3.395){13.69}{60.0}{120.0}
\usefont{T1}{ptm}{m}{n}
\rput(12.892344,7.837188){$ Z$}
\usefont{T1}{ptm}{m}{n}
\rput(15.2,13.6){$ J_p(x_1) \rightarrow \cdots \rightarrow J_p(\tilde{x}_K)$}
\psline[linewidth=0.04cm,arrowsize=0.05291667cm 2.0,arrowlength=1.4,arrowinset=0.4]{<-}(12.1,12.1)(13.5,14.075)
\usefont{T1}{ptm}{m}{n}
\rput(12, 11.8){$ \tilde{x}_1$}
\psline[linewidth=0.04cm,arrowsize=0.05291667cm 2.0,arrowlength=1.4,arrowinset=0.4]{<-}(12,9.75)(11.9,11.5)
\usefont{T1}{ptm}{m}{n}
\rput(12.1, 9.411371){$ x_1$}
\usefont{T1}{ptm}{m}{n}
\rput(11.6,10.8){$P_Z$}
\psline[linewidth=0.04cm,arrowsize=0.05291667cm 2.0,arrowlength=1.4,arrowinset=0.4]{<-}(13.9,13.35)(12.2,9.695)
\usefont{T1}{ptm}{m}{n}
\rput(16.8,11.6){$P_Z\circ J_q^*$}
\usefont{T1}{ptm}{m}{n}
\rput(20.432344,9.971372){$ \|F(x) - y \| \le \hat{\eta}$}
\psarc[linewidth=0.04](22.05,9.145){3.95}{120.0}{210.0}
\end{pspicture}
}
\caption{Projected steepest descent iteration}
\label{fig:proj-landweber}
\end{figure}

\begin{theorem}\label{thm:stopping-rule}
Let Assumption~\ref{assumption:forward-operator} hold true. Moreover, assume that the estimate (\ref{eta}) holds for some positive constant $\eta \in (0, (8 \tilde{C})^{-1})$ 
and $z^\dag \in Z$.

Then Algorithm~\ref{algo:1} stops after a finite number $K = K(\eta)$ of iterations with the discrepancy criterion
\[
   r_K = \|F(x_K)-\yd\| \le \hat{\eta},
\]
being satisfied and strict monotonicity of the Bregman distance
\begin{equation} \label{monoton-Bregman}
\Delta_p(x_{k+1}, z^\dag) \le  \Delta_p(x_{k}, z^\dag) + w_k \Delta_p(x_{k}, z^\dag)^{2/p} -v_k ,
\end{equation}
holds with
\[
 w_k \Delta_p(x_{k}, z^\dag)^{2/p} -v_k < 0,
\]
for all $k \le K(\eta) - 1$.
\end{theorem}

\medskip\medskip

\begin{proof}
We use the same abbreviations for $r_k$ and $t_k$ as in
Algorithm~\ref{algo:1}.

We start with a collection of elementary estimates that will be used
frequently afterwards. With the abbreviations defined in (\ref{abbreviation_2}), (\ref{abbreviation_2a}), inequalities (\ref{Bregman-norm-rela1}) and (\ref{stab-Banach}) yield
\begin{equation}\label{help1}
\begin{aligned}
      \frac{\mathfrak{L}}{2}\|x_k-z^\dagger\|^2 \leq & \frac{\mathfrak{L}}{2}\left( \Delta_p(x_k,z^\dagger) \frac{p}{C_p} \right)^{2/p}
      \\
      \leq & \frac{\mathfrak{L}}{2} \left(\frac{C_p}{p}\right)^{-2/p} \mathfrak{C}^2 \|F(x_k) - F(z^\dag)\|^2 \\
      \leq &\tilde{\mathfrak{C}}(r_k + \|F(z^\dagger)-y\|)^2 \\
      \leq & \tilde{\mathfrak{C}} r_k^2 + 2\tilde{\mathfrak{C}}\eta r_k + \tilde{\mathfrak{C}} \eta^2 \\
      = & r_k - u_k - \eta \;.
\end{aligned}
\end{equation}
With the mean value inequality and (\ref{Bregman-norm-rela1}), it follows that
\begin{equation}\label{help1a}
r_k \le \|F(x_k)-F(z^\dag)\| + \eta \le \hat{\mathfrak{L}} \left(\Delta_p(x_k , z^\dag)\frac{p}{C_p}\right)^{1/p} + \eta .
\end{equation}
Using the definition of $\mu_k$ it follows that for $k=0,1,\ldots$,
\begin{equation}
\label{mu_short}
\mu_k  r_k^{p-1} = \hat{t}_k^{-\frac{1}{q-1}} u_k^{\frac{1}{q-1}} r_k^{p^2-p}\,,\quad
\frac{G_q}{q} \mu_k^q t_k^q = \frac{1}{q} \hat{t}_k^{-\frac{1}{q-1}} u_k^p r_k^{p^2-p} \;.
\end{equation}

Now, we start with the main body of the proof: We claim that
\[
\Delta_p(x_m, z^\dag) < \rho , \quad m = 0,1,\dots,K,
\]
which we prove by induction.
Assume the induction hypothesis that
\[
\Delta_p(x_k, z^\dag) < \rho.
\]
Note that (\ref{converge-radius}) gives the base case. With (\ref{help1a}), we have that
\begin{equation}\label{rkbd}
r_k < \hat{\mathfrak{L}} (\rho \frac{p}{C_p})^{1/p} + \eta = \frac{1+\sqrt{1-8\tilde{\mathfrak{C}} \eta}}{2\tilde{\mathfrak{C}}} -\eta.
\end{equation}
Note that we can rewrite
\[
u_k = -\tilde{\mathfrak{C}} \left(r_k - \frac{1-\sqrt{1-8\tilde{\mathfrak{C}} \eta}}{2\tilde{\mathfrak{C}}} +\eta \right)\left(r_k - \frac{1+\sqrt{1-8\tilde{\mathfrak{C}} \eta}}{2\tilde{\mathfrak{C}}} +\eta \right).
\]
Then, (\ref{rkbd}), combined with the fact that
\[
r_k >\hat{\eta} \ge 3 \eta > \frac{1 - \sqrt{1-8\tilde{\mathfrak{C}} \eta}}{2\tilde{\mathfrak{C}}} -\eta
\]
gives the positiveness of $u_k$. Note that this leads to the positiveness of $v_k$ as following
\[
\begin{aligned}
v_k = & \hat{t}_k^{-\frac{1}{q-1}} u_k^{\frac{1}{q-1}}r_k^{p^2-p} (r_k -\eta - \frac{1}{q} u_k)
\\
\ge & \hat{t}_k^{-\frac{1}{q-1}} u_k^{\frac{1}{q-1}}r_k^{p^2-p} (r_k -\eta - u_k)
\\
= &  \tilde{\mathfrak{C}} \hat{t}_k^{-\frac{1}{q-1}} u_k^{\frac{1}{q-1}}r_k^{p^2-p} (r_k +\eta)^2 >0.
\end{aligned}
\]
Using (\ref{eq:BDpq}) and (\ref{Duality-map-def}) we obtain, for the sequence of residues,
\begin{equation}\label{iterative-bound}
\begin{array}{rl}
   & \Delta_p(\tilde{x}_{k+1},z^\dag)
\\[0.2cm]
   = & \dst{\Delta_p(x_{k},z^\dag) + \frac{1}{q} \left(
     \|\tilde{x}_{k+1}\|^p - \|x_{k}\|^p \right)
         - \langle J_p(\tilde{x}_{k+1}) - J_p(x_{k}) , z^\dag \rangle}
\\[0.2cm]
   = & \dst{\Delta_p(x_{k},z^\dag) + \frac{1}{q} \left(
     \|J_p(\tilde{x}_{k+1})\|^q  - \|J_p(x_{k})\|^q \right)
         - \langle J_p(\tilde{x}_{k+1}) - J_p(x_{k}) , z^\dag \rangle} .
\end{array}
\end{equation}
Applying $(\ref{eq:BDpq})$ and (f) of
Theorem~\ref{thm:Bregman-distance} with $x^* = J_p(\tilde{x}_{k+1})$
and $\tilde{x}^* = J_p(x_k)$, we get
\begin{multline*}
~  \frac{1}{q} \left( \|J_p(\tilde{x}_{k+1})\|^q
                        - \|J_p(x_{k})\|^q \right)\\
\le  \frac{G_q}{q} \|J_p(\tilde{x}_{k+1}) - J_p(x_{k})\|^q
        + \langle J_p(\tilde{x}_{k+1}) - J_p(x_{k}) , x_k \rangle .
\end{multline*}
Substituting (\ref{Steepest descent}) and using this inequality in
(\ref{iterative-bound}) yields
\begin{equation}\label{iterative-bound-2}
\begin{aligned}
   & \Delta_p(\tilde{x}_{k+1}, z^\dag) - \Delta_p(x_{k}, z^\dag)\\
  = & \frac{G_q}{q} \|J_p(\tilde{x}_{k+1}) - J_p(x_{k})\|^q + \langle J_p(\tilde{x}_{k+1}) - J_p(x_{k}), x_k-z^\dagger \rangle\\
  = & \mu_k \left(\frac{G_q}{q} \mu_k^{q-1} t_k^q - \langle T_k ,x_{k} - z^\dag \rangle \right)\;.
\end{aligned}
\end{equation}
We estimate the second term in (\ref{iterative-bound-2}).  Using
(\ref{Bregman-norm-rela1}) and the Lipschitz type stability
(\ref{stab-Banach}), and (\ref{eta}), we find that
\begin{equation}
\label{serious}
\begin{aligned}
   & - \langle T_k , x_{k} - z^\dag \rangle \\
 = & - \langle j_p(R_k) , DF(x_{k})(x_{k} - z^\dag) \rangle \\
 = & - \langle j_p (R_k) , R_k \rangle + \langle j_p(R_k),F(z^\dag) - \yd \rangle )\\
   & \hspace*{2.0cm} + \langle j_p(R_k),F(x_{k}) - F(z^\dag) - DF(x_{k})(x_{k} - z^\dag) \rangle\\
\le & - r_k^{p-1} \left(r_k - \eta - \frac{\mathfrak{L}}{2} \|x_{k} - z^\dag\|^2\right).
\end{aligned}
\end{equation}
From (\ref{iterative-bound-2}) and (\ref{serious}), it follows that, for $k=0,1,2,\ldots$,
\begin{equation}
\label{impro_k_help}
\begin{aligned}
    & \Delta_p(\tilde{x}_{k+1}, z^\dag) - \Delta_p(x_{k}, z^\dag)\\
\le & \mu_k r_k^{p-1} \left(\frac{G_q}{q} \frac{\mu_k^{q-1} t_k^q}{r_k^{p-1}} - r_k + \eta
            + \frac{\mathfrak{L}}{2} \|x_{k} - z^\dag\|^2\right)\; ,
\end{aligned}
\end{equation}
and hence, by (\ref{impro_k_help}), (\ref{Bregman-norm-rela1}) and the non-expansiveness of the Bregman projection (\ref{non-expansive}), we arrive at
\begin{equation}
\label{impro}
\begin{aligned}
& \Delta_p(x_{k+1}, z^\dag) - \Delta_p(x_{k}, z^\dag)\\
\leq  & \Delta_p(\tilde{x}_{k+1}, z^\dag) - \Delta_p(x_{k}, z^\dag)\\
\leq & \mu_k r_k^{p-1} \left(\frac{G_q}{q} \frac{\mu_k^{q-1} t_k^q}{r_k^{p-1}} - r_k + \eta
            + \frac{\mathfrak{L}}{2} \left(\Delta_p(x_k,z^\dag)\frac{p}{C_p}\right)^{2/p}\right).
\end{aligned}
\end{equation}
Using the identities in (\ref{mu_short}) and abbreviations (\ref{abbreviation_2}), (\ref{abbreviation_2a}) we derive that
\begin{equation}
\begin{aligned}
 & \Delta_p(x_{k+1}, z^\dag) - \Delta_p(x_{k}, z^\dag)
\\
\le &  \frac{1}{q} \hat{t}_k^{-\frac{1}{q-1}}u_k^pr_k^{p^2-p} - \hat{t}_k^{-\frac{1}{q-1}} u_k^{\frac{1}{q-1}}r_k^{p^2-p}(r_k - \eta)
\\
& + \frac{\mathfrak{L}}{2}\left(\frac{C_p}{p}\right)^{-2/p}  \hat{t}_k^{-\frac{1}{q-1}}u_k^{\frac{1}{q-1}} r_k^{p^2-p}\Delta_p(x_{k}, z^\dag)^{2/p}
\\
= & -v_k + w_k \Delta_p(x_{k}, z^\dag)^{2/p}.
\end{aligned}
\end{equation}
We finish the proof of the monotonicity of $\Delta_p(x_{k}, z^\dag)$ by showing that
\[
-v_k + w_k \Delta_p(x_{k}, z^\dag)^{2/p} < 0.
\]
In fact,
\begin{equation}
\begin{aligned}
& w_k \Delta_p(x_{k}, z^\dag)^{2/p}
\\
\leq & w_k \mathfrak{C}^2 \|F(x_k) - F(z^\dag)\|^2
\\
\leq &  \frac{\mathfrak{L}}{2}\left(\frac{C_p}{p}\right)^{-2/p} \mathfrak{C}^2 \, (r_k+\eta)^2  \hat{t}_k^{-\frac{1}{q-1}}u_k^{\frac{1}{q-1}} r_k^{p^2-p}
\\
= & (-u_k + r_k - \eta)\hat{t}_k^{-\frac{1}{q-1}}u_k^{\frac{1}{q-1}} r_k^{p^2-p}.
\end{aligned}
\end{equation}
Hence
\begin{equation}\label{help-stop}
\begin{aligned}
& -v_k + w_k \Delta_p(x_{k}, z^\dag)^{2/p}
\\
\le & -v_k - \hat{t}_k^{-\frac{1}{q-1}}u_k^{\frac{1}{q-1} + 1} r_k^{p^2-p} + (r_k - \eta)\hat{t}_k^{-\frac{1}{q-1}}u_k^{\frac{1}{q-1}} r_k^{p^2-p}
\\
= & -\frac{1}{p}\hat{t}_k^{-\frac{1}{q-1}}u_k^{p} r_k^{p^2-p} <0.
\end{aligned}
\end{equation}
The above monotonicity of $\Delta_p(x_{k}, z^\dag)$ with the induction hypothesis completes the induction.

It is left to show that Algorithm~\ref{algo:1} stops after a finite number of iterations(i.e. $K(\eta)$ iterations). We
prove this by contradiction. Suppose that Algorithm~\ref{algo:1} does
not stop within a finite number of iterations and, hence,
\begin{equation}\label{nonstop-rk}
   r_k > \hat{\eta} ,\quad \forall k \ge 0 .
\end{equation}
Then, from the monotonicity of the Bregman distances
(\ref{monoton-Bregman}) and (\ref{help-stop}), we have that
\[
   0 \le \Delta_p(x_{k},z^\dag) \le \Delta_p(x_{0},z^\dag)
              - \frac{1}{p}\sum_{n = 0}^{k - 1} \hat{t}_n^{-\frac{1}{q-1}}u_n^{p} r_n^{p^2-p} ,\quad \forall k > 0 .
\]
It follows that
\[
\sum_{n = 0}^{\infty} \hat{t}_n^{-\frac{1}{q-1}}u_n^{p} r_n^{p^2-p} <\infty
\]
and hence
that $u_k$ converges to $0$ as $k$ goes to infinity. By writing
\[
u_k = -\tilde{\mathfrak{C}} \left(r_k - \frac{1-\sqrt{1-8\tilde{\mathfrak{C}} \eta}}{2\tilde{\mathfrak{C}}} +\eta \right)\left(r_k - \frac{1+\sqrt{1-8\tilde{\mathfrak{C}} \eta}}{2\tilde{\mathfrak{C}}} +\eta \right)
\]
we have that
\[
\lim_{k\rightarrow \infty} r_k = \frac{1-\sqrt{1-8\tilde{\mathfrak{C}} \eta}}{2\tilde{\mathfrak{C}}} - \eta <3\eta <\hat{\eta},
\]
which is a contradiction.

\end{proof}

\medskip\medskip

\begin{remark}
We refer to Algorithm \ref{algo:1} as a steepest descent algorithm in
the sense that it is a generalization of the steepest descent
algorithm for linear inverse problems. Indeed, let $F$ be linear and
assume that we have an unconstrained problem. Then both $\mathfrak{L}$ and $\eta$
can be chosen to be equal to zero. Then we have
\begin{equation*}
   \mu_k = \left( \frac{r_k^p}{t_k^q G_q} \right)^{1/(q-1)} ,
             \quad k=0,1,2,\ldots ,
\end{equation*}
with $r_k = \|F x_k-\yd\|$ and $t_k = \|F^*j_p(F x_k-\yd)\|$. In
particular, for a Hilbert space setting, where
\[
   p=q=2 ,\quad C_p=G_q=1\,, J_p=J_q=Id ,
\]
we get
\begin{equation*}
   \mu_k = \frac{r_k^2}{t_k^2}\,,\quad k=0,1,2,\ldots ,
\end{equation*}
which is the standard parameter choice of the steepest descent method
\cite{Gilyazov1977}. See also \cite{Hanke1991} for efficient adaptations of the Landweber iteration.

In the Hilbert space setting, moreover, the condition (\ref{converge-radius}) requires that $\tilde{\mathfrak{C}} < \frac{1}{8 \eta}$, which in some sense restricts the \emph{curvature}.
Note that for $p=2$ we have $\frac{1}{\|F'(x)\|} \approx \|x-\tilde{x}\|^2/\|F(x)-F(\tilde{x}\|^2 \leq \tilde{\mathfrak{C}}$ and therefore $\frac{\|F''(x)\|}{\|F'(x)\|} \leq \tilde{\mathfrak{C}}{\mathfrak{L}}$, where $\|F'(x)\|$ denotes the operator norm of a directional derivative in direction $x-\tilde{x}$, and $F''$ is the second derivative in the same direction. Thus condition (\ref{converge-radius}) can be interpreted as a curvature to size condition (see \cite{Chavent1996} for the curvature to size concept for variational regularization).
\end{remark}

\medskip\medskip

\begin{remark}
   We refer to (\ref{converge-radius}) as a generalized radius of convergence from the nonlinear Landweber iteration to a steepest descent algorithm in Banach spaces. Indeed, let $\eta$ be equal to zero. Then (\ref{converge-radius}) can be reduced to
  \[
\Delta_p(x_0, z^\dag) < \rho = \hat{\mathfrak{L}}^{-p} \frac{C_p}{p} \tilde{\mathfrak{C}}^{-p} = \left(\frac{C_p}{p}\right)^{3} \left(\cfrac{\hat{\mathfrak{L}} \mathfrak{L} \mathfrak{C}^2}{2}\right)^{-p},
  \]
  which coincides the convergence radius for the nonlinear Landweber iteration in Banach spaces\cite{Hoop2012}.
\end{remark}

\medskip\medskip

\section{Extension to a multi-level algorithm}
\label{sec:progressive}

In this section, we consider a set, $\{ Z_\alpha \}_{\alpha \ge 0}$, of closed and convex subsets of $X$, and an operator
family $\{ F_{\alpha} \}_{\alpha \ge 0}$, where $F_{\alpha}$
is obtained as $F_\alpha = F |_{Z_\alpha}$, or approximations of $F$. We let
\[
\mathcal{B}=\mathcal{B}^{\Delta}_{\rho_0}(x^{\dagger}) = \{ x \in
X\ |\ \Delta_p(x,x^{\dagger}) \le \rho_0 \} \subset \mathcal{D}(F)
 \]
 for some $\rho_0>0$, which is specified in Theorem~\ref{thm:progressive}
and invoke

\medskip\medskip

\begin{assumption}\label{assumption:F-prog}
\begin{enumerate}[(a)]
\item $F$ is weakly sequentially closed, that is,
\begin{equation*}
   \left. \begin{array}{rl}
      & x_n \rightharpoonup x ,
\\
      & F(x_n) \rightarrow y
   \end{array} \right\} \Rightarrow
   \left\{ \begin{array}{rl}
      & x \in \mathcal{D}(F) ,
\\
      & F(x) = y .
   \end{array} \right.
\end{equation*}
\item The Fr\'{e}chet derative, $DF_\alpha$, of $F_\alpha$ is
  Lipschitz continuous on $\mathcal{B} \cap Z_\alpha$ and
\begin{equation}
   \| DF_\alpha(x) \| \le \hat{\mathfrak{L}}_\alpha
   \quad \forall x\in \mathcal{B} \cap Z_\alpha ,
\end{equation}
  \begin{equation}
   \| DF_\alpha(x) - DF_\alpha(\tilde{x}) \|
                  \le \mathfrak{L}_\alpha \| x - \tilde{x} \|
   \quad \forall x, \tilde{x} \in \mathcal{B} \cap Z_\alpha .
\end{equation}
\item The inversion has the uniform Lipschitz type stability for
  elements in $Z_\alpha$, that is, there exists a constant $\mathfrak{C}_\alpha >
  0$ such that
\begin{equation}\label{stab-prog}
   \Delta_p(x,\tilde{x})
        \le \mathfrak{C}_\alpha^p \|F_\alpha(x) - F_\alpha(\tilde{x})\|^p
   \quad \forall x, \tilde{x} \in \mathcal{B} \cap Z_\alpha.
\end{equation}
\end{enumerate}
\end{assumption}

\medskip\medskip

For the stability constants, $\{ \mathfrak{C}_\alpha \}$, and the approximation error, $\{\eta_\alpha\}$, we introduce

\medskip\medskip

\begin{assumption}\label{assumption:index}
\begin{enumerate}[(a)]
\item Let $\eta_\alpha = \eta_\alpha(y)$ be defined by
\[
   \eta_\alpha
        = \operatorname{dist}(y,F_\alpha(Z_\alpha)) ,\quad y \in Y ;
\]
Moreover, we assume that $\eta_\alpha$ is non-negative and monotonically decreasing with respect to $\alpha$ for every fixed
$y \in Y$.
\item If $Z_{\alpha_1} \subset Z_{\alpha_2}$ then $\mathfrak{C}_{\alpha_1} \le
  \mathfrak{C}_{\alpha_2}$.
\item If $\alpha_1 < \alpha_2$ then $Z_{\alpha_1} \subset Z_{\alpha_2}$ and therefore also $\eta_{\alpha_1} \ge \eta_{\alpha_2}$.
\end{enumerate}
\end{assumption}

\medskip\medskip

Typically, the subsets $Z_{\alpha}$ are finite dimensional and the stability constant for the inversion grows with the dimension of these subsets. The
nature of our multi-level algorithm is intimately connected to finding
sparse, albeit approximate, representations of the solution to the
inverse problem, mitigating the mentioned growth of the stability
constants. Indeed, the objective is very similar to \emph{multi-level}
techniques for solving inverse problems \cite{Scherzer1998a,Kaltenbacher2006,Kaltenbacher2008a},
where one exploits that the finite-dimensional problems are stable and
that the outcome of an iteration on a coarse level gives a good
initial guess on a finer level. In this section, we combine any known
controllable factors to an abstract index $\alpha$ of the operator
family and design a progressive iteration method with the aid of the
result from the previous section.

\medskip\medskip

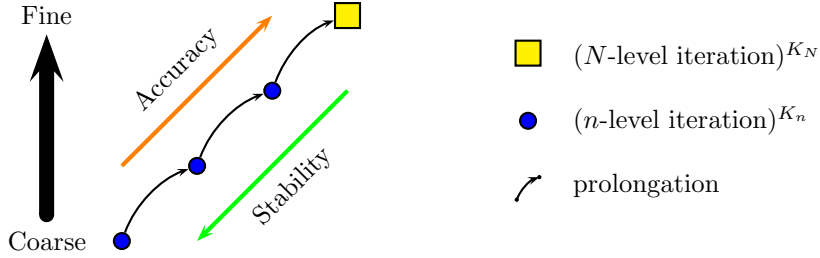
\begin{figure}[ht]
\begin{center}
\unitlength 1cm
\begin{pspicture}(8.8,3.5)

\rput(0,0){Coarse}
\rput(0,3){Fine}
\psline[linewidth=.2cm]{cc->}(0,0.25)(0,2.75)

\cnode[fillstyle=solid,fillcolor=blue](1,0){0.12cm}{P1}
\cnode[fillstyle=solid,fillcolor=blue](2,1){0.12cm}{P2}
\cnode[fillstyle=solid,fillcolor=blue](3,2){0.12cm}{P3}
\fnode[fillstyle=solid,fillcolor=yellow,framesize=0.36cm](4,3){P4}

\pnode(1,1){A1}{}
\pnode(3,3){A2}{}
\ncline[linewidth=1.6pt, linecolor=orange]{->}{A1}{A2}\naput[nrot=:U]{Accuracy}
\pnode(2,0){B1}{}
\pnode(4,2){B2}{}
\ncline[linewidth=1.6pt, linecolor=green]{<-}{B1}{B2}\nbput[nrot=:U]{Stability}

\ncarc[arcangle=25]{->}{P1}{P2}
\ncarc[arcangle=25]{->}{P2}{P3}
\ncarc[arcangle=25]{->}{P3}{P4}

\cnode(6.25,0.55){0.0cm}{A3}
\cnode(6.55,0.85){0.0cm}{A4}
\ncarc[arcangle=25]{->}{A3}{A4}
\put(7.0,0.65){prolongation}


\cnode[fillstyle=solid,fillcolor=blue](6.4,1.6){0.12cm}{S1}
\put(7.0,1.5){$(n$-level iteration$)^{K_n}$}
\fnode[fillstyle=solid,fillcolor=yellow,framesize=0.36cm](6.4,2.5){S2}
\put(7.0,2.35){$(N$-level iteration$)^{K_N}$}
\end{pspicture}
\end{center}
\caption{A illustration of Algorithm~\ref{algo:prog}}
\label{fig:alog}
\end{figure}

In the following algorithm, we refer to the parameter $\alpha$ as an index and only nonnegative integer valued $\alpha$ is considered.

\begin{algorithm}\label{algo:prog}
\begin{enumerate}[$(S1)$]
\setcounter{enumi}{-1}
\item Use $x_{0,0}$ as the starting point. Set $n=0$.
\\[0.1cm]
\item Iteration. Use $F_n$ and $Z_n$ as the modelling operator and
  convex subset to run Algorithm~\ref{algo:1} with the discrepancy
  criterion given by
  \begin{equation}\label{discrepancy-principle-prog}
     K_n = \min \{ k \in\mathbb{N}  \mid
     \| F_n(x_{n,k}) - \yd \| \le (3+\varepsilon)\eta_n \} ,
  \end{equation}
  where $\varepsilon > 0$ is a given uniform tolerance constant.

  STOP, if $n = N$, a given number.
\\[0.1cm]
\item Set $x_{n+1,0} = x_{n,K_n}$, $n = n+1$ and go to
  step $(S1)$.
\end{enumerate}
\end{algorithm}

\medskip\medskip

This algorithm is illustrated in Figure~\ref{fig:alog}.

\medskip\medskip

\begin{theorem}\label{thm:progressive}
Assume that Assumptions~\ref{assumption:F-prog}
and \ref{assumption:index} hold. Assume that there exists a finite
subset of operators, $\{ F_n \}_{n = 1}^N$ say, from the operator
family $\{ F_\alpha \}$ such that
\begin{enumerate}[(a)]
\item The starting point $x_{0,0}$ is within the first convergence
  radius, that is,
\begin{equation}\label{conv-radius-prog}
  \Delta_p(x_{0,0}, z_0^\dag) < \rho_0 ,
\end{equation}
where $z_0^\dag$ denotes the $Z_0$ best approximating solution, i.e.,
\[
   \| F_0(z_0^\dag) - \yd \| = \operatorname{dist}(\yd , F_0(Z_0)) ,
\]
and the $Z_0$ convergence radius $\rho_0$ is defined by
\[
 \rho_0 := \frac{C_p}{p}
\hat{\mathfrak{L}}_0^{-p} \left(\frac{1+\sqrt{1-8\tilde{\mathfrak{C}}_0\eta_0}}{2\tilde{\mathfrak{C}}_0} - 2\eta_0 \right)^p,
\]
with $\tilde{\mathfrak{C}}_{0} = \frac{1}{2}\left(\frac{C_p}{p}\right)^{-2/p} \mathfrak{L}_{0}\mathfrak{C}_{0}^2$;
\\[0.1cm]
\item For every two neighbor levels $Z_n$ and $Z_{n+1}$, $n = 0,\ldots,N-1$, the constants $\eta_n$ and $\eta_{n+1}$, $\hat{\mathfrak{L}}_{n+1}$, $\mathfrak{L}_{n+1}$, $\mathfrak{C}_{n+1}$ satisfy the following inequality
\begin{equation}\label{eta-rho-comp}
   (3 + \varepsilon) \eta_n < \left(\frac{C_p}{p}\right)^{1/p}(\hat{\mathfrak{L}}_{n+1} \mathfrak{C}_{n+1})^{-1} \left(\frac{1+\sqrt{1-8\tilde{\mathfrak{C}}_{n+1}\eta_{n+1}}}{2\tilde{\mathfrak{C}}_{n+1}} - 2\eta_{n+1}\right) -\eta_{n+1} ,
\end{equation}
 where $\tilde{\mathfrak{C}}_{n+1} = \frac{1}{2}\left(\frac{C_p}{p}\right)^{-2/p} \mathfrak{L}_{n+1}\mathfrak{C}_{n+1}^2$.
\\[0.1cm]
\item $N$ is the first positive integer such that $\eta_N \le
  (3+\varepsilon)^{-1} \hat{\eta}$, that is,
\[
   (3 + \varepsilon) \eta_n > \hat{\eta} \quad \forall n < N
\]
and
\[
   (3 + \varepsilon)\eta_N \le \hat{\eta} .
\]
\end{enumerate}
Then, Algorithm~\ref{algo:prog} has the property that it stops after a finite number of iterations when the discrepancy criterion
\begin{equation}\label{discrepancy-final}
   \| F_N(x_{N,K_N}) - \yd \| \le \hat{\eta}
\end{equation}
is satisfied.
\end{theorem}

\medskip\medskip
%
    The strategy of the proof is to estimate the decreasing objective function $\|F(x_{n,K_n}) -\yd\|$ level by level. That
is, one applies Theorem~\ref{thm:stopping-rule} to guarantee that the
discrepancy criterion (\ref{discrepancy-principle-prog}) is attained
with a finite number of iterations on each level $n$. Then, with
(\ref{eta-rho-comp}) and (\ref{discrepancy-principle-prog}), we show
that the initial point $x_{n+1,0}$ on level $n+1$, which
coincides with the iteration result $x_{n,K_n}$ on level $n$,
is within the convergence radius $\rho_{n+1}$. Therefore, the
procedure continues until (\ref{discrepancy-final}) is satisfied.

\begin{proof}
We first adapt the convergence radius, $\rho$, in
Theorem~\ref{thm:stopping-rule} to a $n$-level convergence radius $\rho_n$. For any $n$-level, $n=0,1,2,\dots, N$, one can use Algorithm~\ref{algo:1} to obtain an approximate solution to the operator equation
\[
F_n(x) = y, \quad x\in Z_n,
\]
with a given starting point $x_{n,0}$ and the discrepancy criterion given in (\ref{discrepancy-principle-prog}). If the starting point $x_{n,0}$ satisfy
\begin{equation}\label{n-conv-r}
  \Delta(x_{n,0}, z_n^\dag) < \rho_n := \frac{C_p}{p}\hat{\mathfrak{L}}_{n}^{-p} \left(\frac{1+\sqrt{1-8\tilde{\mathfrak{C}}_{n}\eta_{n}}}{2\tilde{\mathfrak{C}}_{n}} - 2\eta_{n}\right)^p,
\end{equation}
where $z_n^\dag$ denotes the best $Z_n$-approximation, then Theorem~\ref{thm:stopping-rule} can be applied to show that Algorithm~\ref{algo:1} stops after a finite number of iterations with
\[
\| F_n(x_{n,k}) - \yd \| \le (3+\varepsilon)\eta_n
\]
satisfied.
%
Next, we show that, in particular with condition (\ref{eta-rho-comp}), if the starting point for the present level, $x_{n,0}$, is within the convergence radius, then the starting point for the next level, $x_{n+1,0}$, which is equal to $x_{n,K_n}$, is within the convergence radius for the next level. That is to say,
\[
   \Delta_p(x_{n,0}, z_n^\dag) \le \rho_n
\]
implies
\[
\Delta_p(x_{n+1,0}, z_{n+1}^\dag) \le \rho_{n+1},
\]
for all $n < N$.
Indeed, for any $n < N$, according to (\ref{n-conv-r}) and
Theorem~\ref{thm:stopping-rule}, after $K_n$ steps, the $n$-level
discrepancy criterion,
\[
   \| F(x_{n, K_n}) - \yd \|
                    \le (3+\varepsilon) \eta_n ,
\]
is satisfied. Then, with the above inequality and (\ref{stab-prog}), we estimate
\begin{equation}\label{con1}
\begin{aligned}
   & \Delta_p(x_{n+1,0}\, , \, z_{n+1}^\dag)^{1/p}
\\[0.1cm]
   \leq & \mathfrak{C}_{n+1} \|F_{n+1}(x_{n+1,0}) - F_{n+1}(z_{n+1}^\dag)\|
\\[0.1cm]
   \le & \mathfrak{C}_{n+1} (\|F_{n+1}(x_{n+1,0}) - y\| +\|F_{n+1}(z_{n+1}^\dag) - y\|)
\\[0.1cm]
   \le & \mathfrak{C}_{n+1} ( (3+\varepsilon) \eta_n +\eta_{n+1}).
\end{aligned}
\end{equation}
Note that (\ref{eta-rho-comp}) leads to the inequality
\[
\mathfrak{C}_{n+1}((3+\varepsilon) \eta_n +\eta_{n+1}) \le \left(\frac{C_p}{p}\right)^{1/p}\hat{\mathfrak{L}}_{n+1} ^{-1} \left(\frac{1+\sqrt{1-8\tilde{\mathfrak{C}}_{n+1}\eta_{n+1}}}{2\tilde{\mathfrak{C}}_{n+1}} - 2\eta_{n+1}\right)
\]
Substituting this into (\ref{con1}), we have that
\begin{multline}
\Delta_p(x_{n+1,0}\, , \, z_{n+1}^\dag)^{1/p}
\\[0.1cm]
  \leq \left(\frac{C_p}{p}\right)^{1/p}\hat{\mathfrak{L}}_{n+1} ^{-1} \left(\frac{1+\sqrt{1-8\tilde{\mathfrak{C}}_{n+1}\eta_{n+1}}}{2\tilde{\mathfrak{C}}_{n+1}} - 2\eta_{n+1}\right)
=\rho_{n+1}^{1/p}.
\end{multline}
For the last $N$-level, we apply Theorem~\ref{thm:stopping-rule}
again to find that
\[
   \| F_n(x_{N,K_N}) - \yd \|
            \le (3+\varepsilon)\eta_N  \le \hat{\eta} .
\]
\end{proof}

\begin{remark}
  We interpret that Algorithm~\ref{algo:prog} is designed to achieve the optimal (or nearly optimal) accuracy for a feasible starting point. Usually, the finest level bears both the smallest approximation error, which corresponds to the optimal accuracy, and the largest stability constant. Note that the definition of the convergence radius (\ref{converge-radius}) shows its algebraically decaying property with respect to the stability constant. There are cases when only a rough starting point is available. For these cases, one may fail to obtain a reasonable result using Algorithm~\ref{algo:1} directly on the finest level but Algorithm~\ref{algo:prog} leads to a good approximation solution. The condition (\ref{eta-rho-comp}) can be interpreted as a strategy for picking next finer level, which is characterized by its stability constant constants $\mathfrak{C}_{n+1}$, approximation error $\eta_{n+1}$ and $\hat{\mathfrak{L}}_{n+1}$, $\mathfrak{L}_{n+1}$.
\end{remark}

Theorem~\ref{thm:progressive}, especially $(ii)$, indicates that a
sufficient condition for the existence of such a selection of
operators is that the tolerated best-$Z_n$-approximation is within the
convergence radius of $Z_{n + 1}$. In fact, this condition comes from
a bootstrap type competition between $\eta_n$ and $\rho_n$.

We give an example of how conditions (ii) and (iii) in
Theorem~\ref{thm:progressive} can be satisfied.

\medskip\medskip

\begin{example}\label{ex:selection}
Assume that $X$ and $Y$ are Banach spaces and that we can reindex the convex subsets $\{Z_\alpha\}$ such that Assumptions~\ref{assumption:F-prog} and \ref{assumption:index} hold. Moreover, for given $\hat{\eta} > 0$, the following conditions hold:
\begin{enumerate}[(i)]
\item Given starting point $x_{0,0}$ is within the first convergence
  radius $\rho_0$, i.e.,
  \[
  \Delta_p(x_{0,0}, z_0^\dag) < \rho_0 :=  \frac{C_p}{p}
\hat{\mathfrak{L}}_0^{-p} \left(\frac{1+\sqrt{1-8\tilde{\mathfrak{C}}_0\eta_0}}{2\tilde{\mathfrak{C}}_0} - 2\eta_0 \right)^p.
  \]
\item The approximation error $\eta_\alpha = \lambda e^{-\alpha} (\alpha + 2)^{-1}$ for some constant
  $\lambda >> 2\hat{\eta}$.
\item The stability constant $\mathfrak{C}_\alpha = 2 e^\alpha$,
\item The dynamic models of the constants $\hat{\mathfrak{L}}_\alpha$ and $\mathfrak{L}_\alpha$, which are related to the Lipschitz continuity of the Fr\'echet derivative $DF_\alpha$, are given by
    \[
    \hat{\mathfrak{L}}_\alpha = (\alpha + 1)e^{-\alpha} \quad \mbox{and} \quad \mathfrak{L}_\alpha = \tau e^{-\alpha} ,
    \]
for some constant $\tau$ such that
\[
0 < \tau < \left(\frac{C_p}{p}\right)^{3/p} \frac{1}{16 \lambda (4e + 1)}.
\]
\end{enumerate}
Now, we can choose the operators $\{F_n\}_{n=0}^N$ defined by $F_n = F\mid_{Z_n}$ and set the uniform tolerance constant $\varepsilon = 1$ to run Algorithm~\ref{algo:prog}, where $N$ is the first integer such that $
4\eta_N \le \hat{\eta}$ is satisfied. Applying Theorem~\ref{thm:progressive}, we conclude that
\[
   \| F_N(x_{N,K_N}) - \yd \|
                \le \hat{\eta}
\]
is satisfied after a finite number of iterations.

\end{example}
In this example, we can quantify the intermediate constant $\tilde{C}_n$ and the convergence radius $\rho_n$ by
\[
\tilde{\mathfrak{C}}_n = 2\tau  \left(\frac{C_p}{p}\right)^{-2/p} e^{n}
\]
and
\[
\rho_n = \frac{C_p}{p}\hat{\mathfrak{L}}_{n}^{-p} \left(\frac{1+\sqrt{1-8\tilde{\mathfrak{C}}_{n}\eta_{n}}}{2\tilde{\mathfrak{C}}_{n}} - 2\eta_{n}\right)^p.
\]
Noting that
\[
\frac{1}{2} < 1 - 4 \tilde{\mathfrak{C}}_{n}\eta_{n} < 1 + \sqrt{1-8\tilde{\mathfrak{C}}_{n}\eta_{n}} - 4\tilde{\mathfrak{C}}_{n}\eta_{n} < 2 -4\tilde{\mathfrak{C}}_{n}\eta_{n} <2 ,
\]
for $n = 0,1,\dots, N$, we conclude that, for the convergence radius $\rho_n$, the dynamic model is
\[
\left(\frac{C_p}{p} \right)^3 (8 \tau)^{-p} (n+1)^{-p}  < \rho_n  <   \left(\frac{C_p}{p} \right)^3 (2 \tau)^{-p} (n+1)^{-p}.
\]
Let us assume that we are in a situation where only a rough starting point $\tilde{x}$ is available such that
\begin{equation}\label{rs1}
\Delta_p(\tilde{x}, z_0^\dag) < \left(\frac{C_p}{p} \right)^3 (8 \tau)^{-p} < \rho_0
\end{equation}
but
\begin{equation}\label{rs2}
\Delta(\tilde{x}, z_N^\dag) > \left(\frac{C_p}{p} \right)^3 (2 \tau)^{-p} (N+1)^{-p} > \rho_N.
\end{equation}
If we run Algorithm~\ref{algo:1} for single $0$-level, by (\ref{rs1}), Theorem~\ref{thm:stopping-rule} can be applied but the optimal residue estimate we can expect can not be smaller than the $0$-level approximation error $\eta_0 = \lambda /2 >> \hat{\eta}$. If we run Algorithm~\ref{algo:1} for single $N$-level, according to (\ref{rs2}), there is no guarantee that Algorithm~\ref{algo:1} will stop after a finite number of iterations nor yield a reasonable result. Hence a multilevel approach, as Algorithm~\ref{algo:prog}, is proposed to obtain a high-accuracy arroximation $x_{N, K_N}$ satisfying
\[
   \| F(x_{N,K_N}) - \yd \|
                \le \hat{\eta}.
\]

\section{Discussion}
\label{sec:discussion}

We discuss a steepest descent iteration method for solving nonlinear operator
equations in Banach spaces. Provided that the nonlinearity of the forward operator obeys a Lipschitz type stability in a convex and closed subset of the preimage space, we could prove a restricted convergence result and provide an estimate of the error decease. Based on the analysis of the radius of convergence, we introduce a multilevel method and obtain a sufficient condition on the choices of the parameters, mainly on the approximation errors and stability constants.

As an example, we mention inverse boundary value problems for the
Helmholtz equation. Indeed, stability estimates satisfying
Assumption~\ref{assumption:index} have been obtained \cite{Beretta2012} where $Z_{\alpha}$ represents a space spanned by a finite linear
combination of piecewise constant functions. Using our multi-level
algorithm, we arrive at a convergence result by successive
approximation. This result can be further improved, using the same
algorithm, by combining different frequencies and exploiting the
frequency dependence of the stability constants.
The idea of using multiple frequencies was proposed by Chen \cite{Chen1997}, who introduced an algorithm based on recursive linearization. The algorithm
starts with an initial guess at the lowest frequency, which typically
captures the coarse scale variations in the wavepeed. Then the Born
approximation is invoked \cite{Chen1997,Bao2005,Bao2009,Bao2010}. The iteration is based on a linearization of the inverse problem at the present frequency. By progressively increasing the frequency and carrying out the iterations, increasingly finer details are added to the wavespeed model until a sufficiently accurate result is obtained. In \cite{Bao2010}, the convergence of this algorithm was established under certain conditions.
As a direct application of Theorem~\ref{thm:progressive}, the convergence of this algorithm can be revisited. Especially, (\ref{eta-rho-comp}) offers a strategy for picking the frequencies and regularization parameters.

\section*{Acknowledgments}

The research was initiated at the Isaac Newton Institute for
Mathematical Sciences (Cambridge, England) during a programme on
Inverse Problems in Fall 2011.


\bibliographystyle{acm}

\bibliography{conv}%

\end{document}